\documentclass[12pt,  leqno]{article}
\usepackage{amsmath}
\usepackage{amsfonts}
\usepackage{amssymb}
\usepackage{graphicx}
\usepackage{color}
\usepackage{amsthm}
\usepackage{hyperref}

\theoremstyle{plain}

\def\DHLhksqrt#1#2{\setbox0=\hbox{$#1\sqrt{#2\,}$}\dimen0=\ht0
\advance\dimen0-0.2\ht0
\setbox2=\hbox{\vrule height\ht0 depth -\dimen0}%
{\box0\lower0.4pt\box2}}
\newcommand{\bigmu}{\raisebox{-.0\baselineskip}{\Large\ensuremath{\mu}}}

\def\bee{\begin{equation}}
\def\eee{\end{equation}}

\def\DHLhksqrt#1#2{\setbox0=\hbox{$#1\sqrt{#2\,}$}\dimen0=\ht0
\advance\dimen0-0.2\ht0
\setbox2=\hbox{\vrule height\ht0 depth -\dimen0}%
{\box0\lower0.4pt\box2}}

\allowdisplaybreaks

\begin{document}

\bigskip
\bigskip
\bigskip
\bigskip
\bigskip
\bigskip
\centerline{\Large\bf The Disproof of the Riemann Hypothesis}

\bigskip
\bigskip

\begin{center}
{\large \sl{  \sl{ C. Dumitrescu}$^1$, \sl{M. Wolf }$^2$} }  
\end{center}

\begin{center}
$^1$Kitchener,  Canada, email:  cristiand43@gmail.com\\

\medskip

$^2$Cardinal  Stefan  Wyszynski  University, \\  Faculty  of Mathematics and Natural Sciences. \\  
ul. W{\'o}ycickiego 1/3,   PL-01-938 Warsaw,   Poland,e-mail:  primes7@o2.pl
\end{center}

\bigskip
\bigskip

\begin{center}
{\bf Abstract}\\
\bigskip
\begin{minipage}{12.8cm}
We show that there is a contradiction between the Riemann's Hypothesis and some form of the theorem on the  universality of the zeta function.
\end{minipage}
\end{center}

\bigskip\bigskip

\section{Introduction}

In his only paper devoted to the number theory published in 1859 \cite{Riemann-1859} 
(it was also included as an appendix in \cite{Edwards})  Bernhard Riemann  continued  analytically  the series
\bee
\sum_{n=1}^\infty \frac{1}{n^s}, ~~~~~~s=\sigma+it,   ~~~\sigma>1
\eee
to the complex plane  with exception of $s=1$, where the above series is a   divergent harmonic series. He has done it using the integral
\bee
\zeta(s)=\frac{\Gamma(1-s)}{2\pi i}\int_\mathcal{C} \frac{(-z)^s}{e^z-1} \frac{dz}{z}, 
\label{zeta-plane}
\eee
\noindent
\noindent where  the contour $\mathcal{C}$ is \\

\centerline{      }

\centerline{      }

\centerline{      }

\centerline{      }

\centerline{      }

\centerline{      }

\begin{picture}(0,0)(0,0)
\thicklines
\put(200,30){\vector(0,1){70}}
\put(140,60){\vector(1,0){140}}
\put(230,70){\mbox{$\mathcal{C}$}}
\put(200,60){\oval(40,40)[l]}
\put(200,65){\oval(30,30)[rt]}
\put(200,55){\oval(30,30)[rb]}

\put(215,55){\vector( 1,0){45}}
\put(260,65){\vector(-1,0){45}}
\end{picture}

\vskip-1cm
\noindent Appearing in  \eqref{zeta-plane}  the  gamma function  $\Gamma(z)$  has many representations, we present
the Weierstrass product:
\bee
\Gamma(z) = \frac{e^{-C z}}{z} \prod_{k=1}^\infty \left(1 + \frac{z}{k}\right)^{-1} e^{z/k}\,.
\label{Gamma}
\eee
Here  $C$ is the Euler--Mascheroni constant
\[
C = \lim_{n \rightarrow \infty}  \left( \sum_{k=1}^n \frac{1}{k} - \log(n) \right) = 0.577216\ldots.
\label{Euler_Mascheroni}
\]

From (\ref{Gamma}) it is seen that $\Gamma(z)$  is defined for all complex numbers  $z$, except $z = -n$ for integer  $n > 0$,
where are the  simple poles of $\Gamma(z)$.  The most popular definition of the gamma function given by the integral
$\Gamma(z) = \int_0^\infty e^{-t} t^{z-1} dt $  is valid only for $\Re[z]>0$.
Recently perhaps  over 100  representations of  $\zeta(s)$ are known,
for review of the integral and  series  representations see \cite{Milgram-2013}.  The function  $\zeta(s)$ has trivial zeros  at
$s=-2n, ~n=1, 2, 3, \ldots$  and nontrivial zeros  in  the critical strip $0<\Re(s)<1$.  In  \cite{Riemann-1859}   Riemann  made the assumption,  now called the
{\it  Riemann Hypothesis} (RH  for short "in the following), that all nontrivial zeros $\rho$ lie on the {\it critical line} $\Re[s]=\frac{1}{2}$:
$\rho=\frac12+i\gamma$. Often the above requirement is augmented  by the demand that all nontrivial zeros are simple.

There are some results on the zero-free regions of $\zeta(s)$.  K.  Ford  proved  \cite{Ford_2002a}, \cite{FORD2002}  that $\zeta(\sigma+ it)\neq 0$
 whenever
 \bee
 \sigma \geq 1-{\frac {1}{57.54(\log {|t|})^{\frac {2}{3}}(\log {\log {|t|}})^{\frac {1}{3}}}}.
 \label{Ford}
 \eee
 The best present bound seems  to belong to  M. Mossinghoff  and  T. Trudgian \cite{Mossinghoff_Trudgian_2015}:  there are no zeros of
 $\zeta(\sigma + it )$   for  $|t|>2 $  and
 \bee
  \sigma>1-{\frac {1}{5.573412\log |t|}}.
 \eee
This is an improvement over  \eqref{Ford}  as the region free of nontrivial zeros  off critical  line is larger.  In  other words  the width
of the region  with possible  zeros violating the RH  decreases with $t$.

Riemann has shown that $\zeta(s)$ fulfills the {\it  functional identity}:
\bee
\pi^{-\frac{s}{2}}\Gamma\left(\frac{s}{2}\right)\zeta(s)
           =\pi^{-\frac{1-s}{2}}\Gamma\left(\frac{1-s}{2}\right)\zeta(1-s),~~~~{\rm  for } ~ s \in\mathbb{C}\setminus\{0,1\}.
\label{functional-zeta}
\eee
The above form of the functional equation is explicitly symmetrical with respect to the line $\Re(s)=1/2$: the change
$s\rightarrow 1 - s$ on both sides of (\ref{functional-zeta}) shows that the values of the combination
of functions $\pi^{-\frac{s}{2}}\Gamma\left(\frac{s}{2}\right)\zeta(s)$  are  the same at points  $s$  and
$1-s$.  Thus it is convenient to introduce the function
\bee
\xi(s)=\frac{1}{2}s(s-1)\Gamma\left(\frac{s}{2}\right)\zeta(s).
\label{def-xi}
\eee
Then  the functional identity takes the simple form:
\bee
\xi(1-s)=\xi(s)
\label{funct-eq}
\eee

\begin{figure*}[ht]
\includegraphics[page=1,width=.95\textwidth]{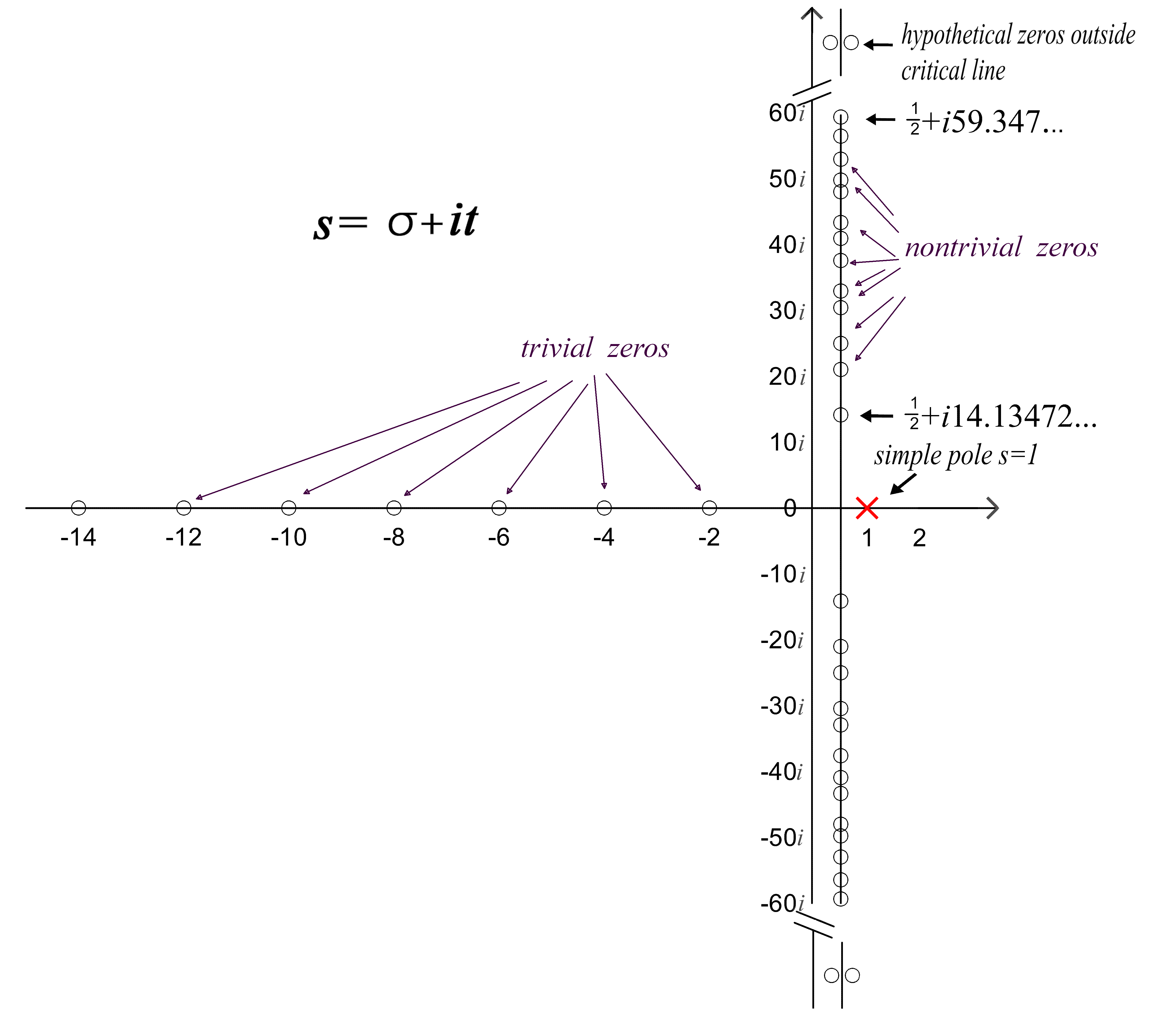}
\caption{\small The location of zeros of the Riemann $\zeta(s)$ function.}
\bigskip
\label{plane}
\end{figure*}

The fact that $\zeta(s)\neq 0$ for $\Re(s)>1$ and the  form  of the  functional
identity entails that {\it nontrivial}  zeros $\rho_n=\beta_n + i\gamma_n$ are located in the {\it critical  strip}:
\[
0 \leq  \Re [\rho_n] =\beta_n \leq 1.
\]
From the complex  conjugation of $\zeta(s)=0$  it follows that if  $\rho_n=\beta_n + i\gamma_n$ is a zero, then
$\overline{\rho_n}=\beta_n -i\gamma_n $ also is a zero.   From the symmetry of the functional equation (\ref{functional-zeta}) with
respect to  the line $\Re[s]=\frac{1}{2}$  it follows, that if $\rho_n=\beta_n + i\gamma_n$ is a zero, then  $1-\rho_n=1-\beta_n - i\gamma_n$  and     $
1-\overline{\rho_n}=1-\beta_n + i\gamma_n$  are also zeros: they are located symmetrically  around  the straight line  $\Re[s]=\frac{1}{2}$ and
the axis $t=0$,  see Fig. \ref{plane}.

Riemann guessed and von Mangoldt proved  that the number of zeta zeros $N(T)$  with positive imaginary parts $<T$ (and real part inside critical
strip)  is  (see eq.(2.3.6)  in  \cite{Borwein_RH})
\bee
N(T)=\frac{T}{2\pi}\log\left(\frac{T}{2\pi}\right)-\frac{T}{2\pi}+\frac78 + \mathcal{O}(\log(T)).
\label{Mangoldt}
\eee

The classical (from XX  century)  references on the RH are \cite{Titchmarsh}, \cite{Edwards},  \cite{Ivic1985}, \cite{Karatsuba-Voronin}.
In the XXI century there appeared two monographs about the zeta function: \cite{Borwein_RH}  and    \cite{Broughan_2017}.\\

There was a  lot of attempts to prove RH and the common  opinion  was that it is  true.  However let us notice that there were  famous
mathematicians: J. E. Littlewood \cite{linki9},  P. Turan and A.M. Turing  \cite[p.1209]{Booker_turingand} believing  that the  RH is not true.
In Karatsuba's talk  \cite{linki_Karatsuba} at 1:01:10 \footnote{We thank A. Kourbatov for bringing  this fact to our  attention}
he mentions that   Atle Selberg had serious doubts whether RH is true or not.
See also  the paper ``On some reasons for doubting the
Riemann  hypothesis'' \cite{Ivic2003}  (reprinted in \cite[p.137]{Borwein_RH})  written by  A. Ivi{\'c}.
New arguments against RH   can be found in \cite{Blanc_2019}.  There were  some attempts to prove RH using the  physical
methods,  see  \cite{Physics-and-RH-RevModPhys, Wolf_2020_RPP}.

Some  analogs of  the RH were proved and some other  were disproved.  Andr{\'e} Weil proved the Riemann hypothesis to be true for field functions
\cite{Weil_1948},  while  the Davenport-Heilbronn zeta-function \cite{Davenport_Heilbronn_1936},  which shares many properties with usual $\zeta(s)$,  has zeros  outside critical line   and even to the right of $\Re(s)=1$,  see \cite{Balanzario_Sanchez-Ortiz_2007}.

In this paper we are going to present  incompatibility  of  RH with the  theorem  on the  universality of  the logarithmic  derivative of
the  zeta function.    In 1975   S.M. Voronin   \cite{Voronin-theorem}  proved first universality theorem.  He wanted to prove the RH but
instead he proved  remarkable theorem on the universality  of the Riemann $\zeta(s)$ function  (see  also \cite[p.251]{Karatsuba-Voronin}):

\bigskip

{\bf Voronin's theorem}:  Let $0<r<\frac14$,  Suppose that the function  $F(s)$  is analytic  for  $|s|<r$  and continous   for  $|s|\leq r$.
If  $F(s)\neq 0$  for  $|s|<r$,  then for  any  $\epsilon > 0$   there exists a $T(\epsilon,F)$  such that
\bee
\max_{|s| \leq r} \; \left| F(s) - \zeta\left(s+\left(\frac{3}{4} + i\;T(\epsilon,f)\right)\right)  \right| < \epsilon.
\label{Voronin}
\eee

Put simply in words it means that the zeta function approximates  locally any smooth function in a uniform way.

We will use  the  universality theorem  for the  logarithmic  derivative  of $\zeta(s)$  \cite[Theorem 2]{Laurincikas1985ZerosOT}:

{\bf Laurin{\v{c}}ikas's theorem}:  Let  $K$  be a compact subset of the strip $D(1/2,1)=\{s\in \mathbb{C}: \frac12<\Re(s)<1\}\subset \mathbb{C}$
such that the complement of $K$  is connected. Let $F : K \to   \mathbb{C}$  be an analytic  function inside $K$ and continuous up to the boundary
of $K$. Then for  any  $\epsilon > 0$:
\bee
\liminf_{T\to\infty}\frac{1}{2T}{\bigmu}\left\{{\tau\in[-T, T]}:\sup_{s\in K}\left|\frac{\zeta'}{\zeta}(s+i\tau)-F(s)\right|<\epsilon\right\}>0.
\label{Litwin}
\eee

Above  $\mu(A) $ is the Lebesgue measure of the set $A \subset \mathbb{R}$.  Let us remark that in Voronin's theorem the existence of just one $T$  such
that  $\zeta(s+iT)$  approximates $F(s)$  is  guaranteed  but in  Laurin{\v{c}}ikas's theorem  the existence  of uncountably  many such $T$  is
secured.  There are also  such ``strong'' or `` enhanced'' versions of the Voronin  theorem,  where the  set of  $T$  is uncountable,  see
\cite{Laurincikas_1996},  \cite{MATSUMOTO_2015}.
In fact,  because we assume  RH,  $\zeta'(s)/\zeta(s)$  is  continuous  to the right of $\Re(s)=\frac12$.   there  will be  some (short)
{\it  intervals}  of  $T$  such that $\zeta'(s+iT)/\zeta(s+iT)$  approximates $F(s)$: if some  $T'$ is good
then also  neighboring  $T$  will satisfy  $ |\zeta'(s+iT)/\zeta(s+iT)-F(s)|<\epsilon$.
We  stress  that  in these  theorems  the {\it  RH is not assumed}.  The idea of our proof is as follows:  We  assume  RH,  thus inside  any
rectangle  on the right of critical line there  should be   no  zeros  of $\zeta(s)$.  Using the  Argument  Principle  and  the universality of
$\zeta'(s)/\zeta(s)$  we will  show  that  it is possible to  choose   such  rectangles that  integral  of  $\zeta'(s)/\zeta(s)$  along  edges  of  these
rectangles will  be  different  from  zero  meaning  that  there  are  zeros of  $\zeta(s)$  to   the  right  of  $\Re(s)=\frac12$:   contradiction.
Thus our  main result  is:

\bigskip

{\bf Theorem}:  The Riemann Hypothesis  is not   true.

\bigskip

In the following section we will  present  the  proof  of above theorem.

\section{Incompatibility between Riemann's Hypothesis and universality of $\zeta'(s)/\zeta(s)$}

We will use the    Argument Principle, see  e.g. \cite[\S 4.10]{Henrici_v1}.  or  [Sect.2.3]\cite{Borwein_RH}:

\bigskip

\noindent {\bf The Argument Principle.} If $F(z)$ is a meromorphic function inside and on some closed contour $\mathcal{C}$, and $F$ has no zeros
or poles on $\mathcal{C}$, then
\bee
{\frac {1}{2\pi i}}\oint _{\mathcal{C}}\frac{F'}{F}(z)\,dz=Z-P
\label{argument_principle}
\eee
where $Z$ and $P$ denote respectively the number of zeros and poles of $F(z)$  inside the contour $\mathcal{C}$, with each zero and pole counted as
many times as its multiplicity.

\bigskip

Let  us  notice that the only pole of $\zeta(s)$  is  at  $s=1$.  We consider the rectangle $ABCD$ shown in Figure 2.  We will denote this rectangle
also as  $\mathcal{D}(\alpha, \beta ,T)$  where  $\frac12<\alpha<\beta<1$  and $T$  will  be determined later.   If we assume that the
RH is true,  then inside  this rectangle  and on its  border  there are no zeros of $\zeta(s)$,
thus   the integral
\bee
\oint _{ABCD}\frac{\zeta'}{\zeta}(s)\,ds
\label{calka1}
\eee
should vanish  (for  $\beta<1$  pole  $s=1$  is outside  $\mathcal{D}(\alpha, \beta ,T)$).  Using the  strong   universality  of $\zeta'(s)/\zeta(s)$
we will show that it is possible to choose such $T$  that this integral will be different from  zero.

\setlength{\unitlength}{1.4cm}
\begin{picture}(10,10)(-5,-4)
\thicklines
\put(-3.5,0){\vector(1,0){6}}
\put(2.7,-0.1){$\sigma$}
\put(-0.9,0.1){$\alpha$}
\put(1.4,0.1){$\beta$}
\put(2.2,0.1){$1$}
\put(2.1,-0.1){\line(0,1){0.2}}
\put(-3.0,-4.5){\vector(0,2){9}}
\put(-1.1,-4.5){\color{red}\line(0,1){8.9}}
\put(-1.0,4.2){${\small \frac12 +it}$}
\put(-0.6,3.2){${\small \alpha +iT}$}
\put(1.35,3.4){${\small \beta +iT}$}
\put(-2.9,0.1){$0$}

\multiput(1.3,-4)(0,2.5){3}{\vector(0,1){2.5}}
\multiput(-0.7,-1.5)(0,2.5){3}{\vector(0,-1){2.5}}
\put(1.3,3.5){\vector(-1,0){2}}
\put(-0.7,-4){\vector(1,0){2}}
\put(-0.59,-3.9){${\small\displaystyle{\small{\alpha-iT}}}$}
\put(1.33,-3.9){${\small\displaystyle{\small{\beta-iT}}}$}

\multiput(2.1,-4)(0,0.8){11} {\line(0,1){0.2}}

\put(1.1,3.6){$A$}
\put(-1,3.5){$B$}
\put(-1,-4.2){$C$}
\put(1.3,-4.2){$D$}

\put(-2.9,4.5){$t$}
\put(-4, -5){Figure 2. ~The  rectangle $\mathcal{D}(\alpha, \beta,  T)$,  in red critical line is plotted.}
\put(-4.8, -5.4) { Sides  $AB$ and $CD$  have length  $\beta-\alpha$  while  sides $CB$ and $DA$  have length $2T$. }\\
\end{picture}

\bigskip
\bigskip

\centerline{      }

\centerline{      }

\centerline{      }

\bigskip
\bigskip

We split  the  integral \eqref{calka1}   into four  integrals
\bee
\oint _{ABCD}\frac{\zeta'}{\zeta}(s)\,ds=\int _{DA}\frac{\zeta'}{\zeta}(s)\,ds+\int _{AB}\frac{\zeta'}{\zeta}(s)\,ds+
\int _{BC}\frac{\zeta'}{\zeta}(s)\,ds +\int _{CD}\frac{\zeta'}{\zeta}(s)\,ds.
\label{calka2}
\eee
We choose  $F(s)$ from   the Laurin{\v{c}}ikas's theorem to be the complex constant  on a compact $K$  on the right half of the critical strip
(that $K$ can be chosen a horizontal segment):
\bee
F(s)=U+iV,  ~~~~~U, V \in \mathbb{R}.
\eee
This  $F(s)$  in fact  can be  constant   on the whole complex plane $\mathbb{C}$ We can put $U=0$  as it will cancels  out  later.  From the
Laurin{\v{c}}ikas's theorem  we  know  that there exists such $T$  that
\bee
\frac{\zeta'}{\zeta}(s+iT)=U+iV +r(s)
\eee
where  the  remainder   $|r(s)|<\epsilon_1$  (later we will use  yet other  epsilons)  on the horizontal side of our rectangle $(\alpha+iT, \beta+iT)$,  i.e. we  choose  $T$  as the height of our rectangle  $ABCD$.
Thus for segment $AB$  we have
\bee
\int _{AB}\frac{\zeta'}{\zeta}(s)\,ds=(\alpha-\beta)(U+iV) + o(1)
\eee
where $o(1)$  absorbs  term:
\bee
\left|\int _{AB} r(s)\,ds\right|<(\beta-\alpha)\epsilon_1
\eee
and  can be made arbitrarily small.  Because there is a symmetry  with respect to real axis:
\bee
\frac{\zeta'}{\zeta}(\overline{s})=\overline{\frac{\zeta'}{\zeta}(s})
\eee
we have
\bee
\int _{CD}\frac{\zeta'}{\zeta}(s)\,ds=(\beta-\alpha)(U-iV) + o(1)
\eee
and together
\bee
\int _{AB}\frac{\zeta'}{\zeta}(s)\,ds + \int _{CD}\frac{\zeta'}{\zeta}(s)\,ds =2i(\alpha-\beta)V + o(1).
\eee

\bigskip

To calculate  integrals along vertical segments we will use the  following   formula for the logarithmic  derivative of the
$\zeta(s)$  function (see e.g. \cite[Chapt.12]{Davenport}):
\bee
\frac{\zeta'(s)}{\zeta(s)}=\frac{1}{1-s}+\frac12 \log(\pi)-\frac12 \frac{\Gamma'(\frac{s}{2}+1)} {\Gamma(\frac{s}{2}+1)}+\sum_\rho\Big(\frac{1}{s-\rho}\Big).
\label{prop1}
\eee

The last summation  extends  over all  nontrivial zeros of $\zeta(s)$,  i.e.  assuming the RH  over $\rho=\frac12\pm i\gamma$.

\bigskip

Integrating  the first  term  on  rhs of \eqref{prop1}  we obtain
\bee
\int_{DA} \frac{1}{1-s} \,ds = \log(\beta-1+iT)-\log(\beta-1-iT)
\eee
We have
\bee
\log(z)=\log|z|+i\arg(z)
\label{logarytm}
\eee
hence we obtain
\bee
\int_{DA} \frac{1}{1-s}\,ds = 2i\arg(\beta-1+iT)
\eee
Similarly for left vertical segment $CB$  we obtain
\bee
\int_{CB} \frac{1}{1-s}\,ds = 2i\arg(\alpha-1+iT)
\eee
Together for  the positive circulation of the contour $ABCD$  we have
\bee
\int_{DA} \frac{1}{1-s}\,ds - \int_{CB} \frac{1}{1-s}\,ds = 2i(\arg(\beta-1+iT)-\arg(\alpha-1+iT))
\eee
The  Laurin{\v{c}}ikas's theorem assures that $T$ can jump to  arbitrarily  large  $T$  hence:
\bee
\lim_{T\to \infty} \arg(\alpha-1+iT)=\frac{\pi}{2}
\eee
\[
\lim_{T\to \infty} \arg(\beta-1+iT)=\frac{\pi}{2}
\]
and  it  follows that  for  sufficiently  large $T$:
\bee
\int_{DA} \frac{1}{1-s}\,ds +\int_{BC} \frac{1}{1-s}\,ds = o(1).
\eee
For the constant term  $\frac12 \log(\pi)$  in  \eqref{prop1}  we obtain:
\bee
\int_{DA}\frac12 \log(\pi)\, ds = \frac12 \log(\pi)\int_{-T}^T i \, dt = \log(\pi)T\, i
\eee
Adding  integral over  left vertical segment we obtain:
\bee
\int_{DA}\frac12 \log(\pi)\, ds + \int_{BC}\frac12 \log(\pi)\, ds =0
\eee

Next we will calculate  integral of the logarithmic derivative of $\Gamma(s)$   function.  We recall that
\bee
\frac{\Gamma'(s)}{\Gamma(s)}=\psi(s),
\eee
where $\psi(s)$ is the digamma function. In \cite[formula  (6.3.18)]{abramowitz+stegun}  we found:
\bee
\psi(s)=\log(s)-\frac{1}{2s}-\frac{1}{12s^2}+\ldots ~~~s\to \infty, ~~|\arg(s)|<\pi.
\label{psi_od_z}
\eee
We note that
\bee
\frac{d}{ds}\log(s+2)=\frac{1}{s+2},
\eee
\bee
\frac{d}{ds}\left((s+2)\log(1+\frac{s}{2})-(s+2)\right)=\log\left(1+\frac{s}{2}\right),
\eee
thus  we  have:
\bee
\int_{DA}  \frac{\Gamma'(\frac{s}{2}+1)}{\Gamma(\frac{s}{2}+1)}\, ds=\left(-(s+2)+(2+s)\log\left(1+\frac{s}{2}\right)-\log\left(2+s\right)\right)\Big|_{\beta-iT}^{\beta+iT}+\mathcal{O}\left(\frac{1}{T}\right).
\label{LogGamma33}
\eee
Because  $1+\alpha/2>0$  and hence also $1+\beta/2>0, ~2+\alpha>0, ~2+\beta>0  $  we avoid branch cut in logarithm for negative arguments.
After straightforward  calculations we obtain for $T\to \infty$  (in fact $T$  does not increase to infinity continuously but
$T$  leaps  over  large  values in according   with   Laurin{\v{c}}ikas's theorem):
\bee
\frac12 \left\{ \int_{DA}  \frac{\Gamma'(\frac{s}{2}+1)}{\Gamma(\frac{s}{2}+1)}\, ds + \int_{BC}  \frac{\Gamma'(\frac{s}{2}+1)}{\Gamma(\frac{s}{2}+1)}\, ds\right\}\xrightarrow{~~T \to \infty~~} (\beta-\alpha)\frac{\pi}{2}i.
\eee

Now we will calculate the integrals of the last sum in \eqref{prop1}  over  vertical edges of the rectangle $ABCD$.   In the sum
\bee
S(t)=\sum_\rho \frac{1}{s-\rho},
\label{sumka}
\eee
where $s\in (\alpha-iT, \alpha+iT)$, i.e. $s=\alpha+it$.  The sum  over zeros $\rho $ is convergent  when zeros  $\rho$ and complex conjugate
$\overline{\rho}$ are  paired together  and under RhH we have
\bee
\frac{1}{s-\frac12-i\gamma_n}+\frac{1}{s-\frac12+i\gamma_n}=2\frac{s-\frac12}{(s-\frac12)^2+\gamma_n^2}.
\eee
In  $S(t)$   the zeros are  summed according to   increasing absolute  values of the imaginary parts $\gamma_n$,  see \cite[Sect.2.5]{Edwards}
(it follows  also   from the M-test of Weierstrass):  the sum  $\sum_\rho 1/|\rho-1/2|^{1+\epsilon}$  over all roots  of $\zeta(s)$  converges for any $\epsilon>0$.  If RH is not true there is a
possibility  that for some zero $\rho'$  denominator $s-\rho'$ will
be  zero and the sum   will be infinite.  Only when we assume  RH $\rho=\frac12+i\gamma$,  thus even if for some $k$  there will be  $t=\gamma_k$
the denominator  will be $\alpha-\frac12\neq 0$  and we  have uniformly in $t$  that  for any  $\epsilon_2>0$  there exists such
$N(\epsilon_2)$  that
\bee
\left|\sum_{n>N(\epsilon_2)} \frac{1}{(s-\frac12)^2+\gamma_n^2}\right|<\epsilon_2 
\eee
for $s$ on the vertical edges of rectangle: $s\in (\alpha-iT, \alpha+iT)$  or $s\in (\beta-iT, \beta+iT)$.  To  have negligible  integrals of the
above  sum over zeros  along the edges $DA$  and  $BC$  we  choose
\bee
\epsilon_2 =\frac{1}{T^2}.
\label{epsylon2}
\eee
Thus  the remaining sum is finite and we can interchange integration with the  summation.
For the vertical  $CB$  we write
\bee
\int_{BC} \sum_{n\leq N(\epsilon_2)} \frac{1}{s-\frac12-i\gamma_n}+\frac{1}{s-\frac12+i\gamma_n}\, ds=
\eee
\[
\sum_{n\leq N(\epsilon_2)} \Bigg(\log\left(\alpha-\frac12 -iT-i\gamma_n\right)-\log\left(\alpha-\frac12+iT-i\gamma_n\right)+
\]
\[
~~~~~~~~~~~~~~~\log\left(\alpha-\frac12 -iT+i\gamma_n\right)-\log\left(\alpha-\frac12+iT+i\gamma_n\right)\Bigg)
\]
The real parts of arguments under  logarithms are positive and  we  use   the  main  branch.   Now we make use of \eqref{logarytm}
and the fact that $\log(|z|)=\log(|\overline{z}|)$  and $\arg(\overline{z})=-\arg(z)$.  Adding  similar
expression  for edge $DA$     we  obtain for the summands:
\bee
\int_{BC}  \frac{1}{s-\frac12-i\gamma_n}+\frac{1}{s-\frac12+i\gamma_n}\, ds + \int_{DA}  \frac{1}{s-\frac12-i\gamma_n}+\frac{1}{s-\frac12+i\gamma_n}\, ds =
\eee
\[
2i \Bigg(\arg\left(\beta-\frac12 +i(T-\gamma_n)\right)+\arg\left(\beta-\frac12+i(T+\gamma_n)\right)
\]
\[
~~~~~~~~~~~~~~-\arg\left(\alpha-\frac12 -i(T-\gamma_n)\right)-\arg\left(\alpha-\frac12+i(T+\gamma_n)\right)\Bigg)
\]
Next we use the relation
\bee
\arg(x+iy)=\arctan\left(\frac{y}{x}\right)~~~\mbox{for ~~} x>0
\eee
and obtain:
\bee
\sum_{n\leq N(\epsilon_2)} \Bigg(\int_{BC}  \frac{1}{s-\frac12 + i\gamma_n}\ +  \frac{1}{s-\frac12 - i\gamma_n}\, ds
\eee
\[
+\int_{DA}  \frac{1}{s-\frac12 + i\gamma_n} + \frac{1}{s-\frac12 - i\gamma_n} \, ds \Bigg) =
\]
\[
2i  \sum_{n\leq N(\epsilon_2)} \Bigg(\arctan\left(\frac{T-\gamma_n}{\beta-\frac12}\right)+\arctan\left(\frac{T+\gamma_n}{\beta-\frac12}\right)
\]
\[
-\arctan\left(\frac{T-\gamma_n}{\alpha-\frac12}\right)-\arctan\left(\frac{T+\gamma_n}{\alpha-\frac12}\right)\Bigg)
\]
We collect  separately  terms  with  $T+\gamma_n$  and  $T-\gamma_n$  and  our sum takes  the  following  form  ($N=N(\epsilon_2)$  for short):
\bee
S_N:=\sum_{k=1}^N \Bigg( \arctan\left(\frac{T-\gamma_k}{\beta-\frac12}\right)-\arctan\left(\frac{T-\gamma_k}{\alpha-\frac12}\right)
\label{suma_1}
\eee
\[
+ \arctan\left(\frac{T+\gamma_k}{\beta-\frac12}\right)-\arctan\left(\frac{T+\gamma_k}{\alpha-\frac12}\right)\Bigg)
\]

Fortunately there is a method of calculating sums of $\arctan$  functions by telescoping.  Namely we have  (see  [eq.(2.5)]\cite{Boros_Moll_2005}):

\bigskip

{\bf Lemma  1}:
Let $f(x)$  be  defined by
\bee
h(x)=\frac{f(x+1)-f(x)}{1+f(x+1)f(x)}
\eee
Then
\bee
\sum_{k=1}^n \arctan(h(k))=\arctan(f(n+1))-\arctan(f(1))+\pi\sum {\rm sgn}~f(k)
\eee
where the last  sum is taken over all  $1\leq k\leq n$ for which $f(k+1)f(k)<-1$   and
where
\[
{\rm sgn}(x)=\begin{cases}
1 & \text{for ~~} x>0  \\
0 & \text{for ~~} x=0  \\
-1 &  \text{for~~} x<0.
    \end{cases}
\]
\bigskip

 We have  \cite[formula (4.4.34)]{abramowitz+stegun}
\bee
\arctan(x)+\arctan(y)=\left\{
\begin{array}{ll}
\displaystyle{\arctan\left(\frac{x+y}{1-xy}\right)} & \mbox {if}~~  xy<1 \\
\medskip\\
\displaystyle{\arctan\left(\frac{x+y}{1-xy}\right)+\pi{\rm sgn}(x) } & \mbox {if~~} xy>1
\end{array}
\right.
\eee

Using it in \eqref{suma_1}  we obtain  terms
\bee
\arctan\left(\frac{T-\gamma_k}{\beta-\frac12}\right)-\arctan\left(\frac{T-\gamma_k}{\alpha-\frac12}\right)=
\eee
\[
\arctan\left(\frac{(\alpha-\beta)(T-\gamma_k)}{(\alpha-\frac12)(\beta-\frac12)+(T-\gamma_k)^2}\right)+\pi{\rm sgn}\left(\frac{T-\gamma_k}{\beta-\frac12}\right),
\]
\bee
\arctan\left(\frac{T+\gamma_k}{\beta-\frac12}\right)-\arctan\left(\frac{T+\gamma_k}{\alpha-\frac12}\right)=
\eee
\[
\arctan\left(\frac{(\alpha-\beta)(T+\gamma_k)}{(\alpha-\frac12)(\beta-\frac12)+(T+\gamma_k)^2}\right)+\pi{\rm sgn}\left(\frac{T+\gamma_k}{\beta-\frac12}\right)
\]
Thus we have the sum
\bee
S_N=\sum_{k=1}^N \arctan\left(\frac{(\alpha-\beta)(T-\gamma_k)}{(\alpha-\frac12)(\beta-\frac12)+(T-\gamma_k)^2}\right) +
\eee
\[
\sum_{k=1}^N \arctan\left(\frac{(\alpha-\beta)(T+\gamma_k)}{(\alpha-\frac12)(\beta-\frac12)+(T+\gamma_k)^2}\right) +Q_1\pi
\]
where $Q$  is an integer.  We introduce two functions:
\bee
h_1(k)=\frac{(\alpha-\beta)(T-\gamma_k)}{(\alpha-\frac12)(\beta-\frac12)+(T-\gamma_k)^2}  
\eee
\bee
h_2(k)=\frac{(\alpha-\beta)(T+\gamma_k)}{(\alpha-\frac12)(\beta-\frac12)+(T+\gamma_k)^2} 
\eee
To telescope sums \eqref{suma_1}  we  need to find two functions $f(k)$ and  $g(k)$   fulfilling  recurrences:
\bee
f(k+1)=\frac{(\alpha-\beta)(T-\gamma_k)+(\alpha-\frac12)(\beta-\frac12)f(k)+(T-\gamma_k)^2f(k)}{(\alpha-\frac12)(\beta-\frac12)+(T-\gamma_k)^2-(\alpha-\beta)(T-\gamma_k)f(k)},,
\label{rek1}
\eee
\bee
g(k+1)=\frac{(\alpha-\beta)(T+\gamma_k)+(\alpha-\frac12)(\beta-\frac12)g(k)+(T+\gamma_k)^2g(k)}{(\alpha-\frac12)(\beta-\frac12)+(T+\gamma_k)^2-(\alpha-\beta)(T+\gamma_k)g(k)}
\label{rek2}
\eee
with the initial conditions $~f(1)=0,  ~g(1)=0$.   We do not need explicit solutions of above equations  just the behavior  of these  functions for large $k$.
In the Appendix we will prove the

\bigskip

{\bf Lemma 2}:  We have
\bee
\lim_{N \to \infty} f(N)=\infty
\eee
\bee
\lim_{N \to \infty} g(N)=-\infty
\eee

Hence:
\bee
\lim_{N \to \infty} \arctan f(N)=\frac{\pi}{2}
\eee
\bee
\lim_{N \to \infty} \arctan g(N)=- \frac{\pi}{2}
\eee
Larger  $N(\epsilon_2)$  means smaller $\epsilon_2$  and smaller term $o(1)$.  Finally we obtain:
\bee
S_N=Q\pi+o(1)
\label{suma61}
\eee
where $Q$ is   an  integer.  For the total  integral  over the  rectangle  $ABCD$  we get:
\bee
\frac{1}{2\pi i}\oint_{ABCD} \frac{\zeta'}{\zeta}(s) ds = \frac{\beta-\alpha}{4} + \frac{\alpha-\beta}{\pi}V+Q+o(1).
\label{calka2}
\eee
We  can choose  $\alpha$  and $\beta$  ($\frac12< \alpha<\beta<1$)    and $V$  as  we wish  and we   put
\bee
\alpha=\frac35, ~~~~\beta=\frac45, ~~~~V=- \pi.
\label{war1}
\eee
and it   imposes  condition
\bee
\frac{\beta-\alpha}{4} + \frac{\alpha-\beta}{\pi}V=\frac14
\eee
We also choose  parameters such that  the modulus of the $o(1)$  term is strictly less than $\frac12$.  At this  point we reach
contradiction:  we assumed RH hence the above  integral \eqref{calka2} should be zero, but we  have  impossible equality
with  integer $Q$:
\bee
\frac14  + Q + o(1) = 0.
\eee
Thus RH cannot be true.

In general  we  can  take
\bee
\alpha=\frac35, ~~~~\beta=\frac45, ~~~~V=(\frac14 -5r) \pi.
\eee
and then we  get
\bee
\frac{\beta-\alpha}{4} + \frac{\alpha-\beta}{\pi}V=r
\eee
where $r$  is arbitrary real   number,  while  from the   Principle  Argument the integral \eqref{calka2}  over $ABCD$  should be an integer,
but we reached  this  {\bf contradiction}  under the {\it  assumption of RH}:  there are no zeros  of $\zeta(s)$  on the rectangle   $ABCD$
and  we were able to calcluate the sum  \eqref{sumka}  only under  the assumption  that all $\rho=\frac12 +i\gamma$.

\section{Final  Comments}

We assumed RH and  using the   Laurin{\v{c}}ikas's universality   theorem we reached the contradiction:
in \eqref{argument_principle}  the integral on lhs  (divided by $2\pi i$)  has to be an  integer  but  under RH  we can  choose  parameters
$\alpha,  ~\beta, ~V$  such that we will  get  arbitrary  real  number,   not an integer.
.  In  our final  expression  \eqref{calka2}  with the parameters we chose, this value is not an integer.

 There are  effective versions  of  the  universality  Voronin's  Theorem   e.g. \cite{Garunkstis_2010},
and  obtained  numbers are  enormously  large:  they  grow like  double iterated  exponential  functions.  We can  expect  that $T$  appearing in  the
Laurin{\v{c}}ikas's theorem is also very large.  Thus we can say   that   RH  is  {\it practically  true}.
 Odlyzko has expressed the view that the hypothetical zeros off
the critical line are unlikely to be encountered for  $t$  below $10^{10^{10000}}$, see \cite[p.358]{Derbyshire}.
It seems  there is no  hope to  find explicit values   of the zeros
off the critical  line.

In \cite{Holt2018}  on pp.60--61 we can read: ``Long open conjectures in analysis tend to be false ''.

\bigskip
\bigskip
{\bf Acknowledgment:}  We  thank  Alexei  Kourbatov for  many  remarks improving the manuscript.
|

\bigskip
\bigskip

\section{Appendix: proof  of Lemma  2}

\bigskip

We repeat here   in  another form  equations  \eqref{rek1},  \eqref{rek2}
\bee
f(k+1)=\frac{\left((T-\gamma_k)^2+(\alpha-\frac12)(\beta-\frac12)\right)f(k)+(\alpha-\beta)(T-\gamma_k)}{-(\alpha-\beta)(T-\gamma_k)f(k)+(T-\gamma_k)^2+(\alpha-\frac12)(\beta-\frac12)},
\label{rek3}
\eee
\bee
g(k+1)=\frac{\left((\alpha-\frac12)(\beta-\frac12)+(T+\gamma_k)^2\right)g(k)+(\alpha-\beta)(T+\gamma_k)}{-(\alpha-\beta)(T+\gamma_k)g(k)+(T+\gamma_k)^2+(\alpha-\frac12)(\beta-\frac12)},
\label{rek4}
\eee
with initial conditions $~f(1)=0, g(1)=0$.  The above equations have the  form of the   Riccati  equation, see e.g  \cite[Appendix A]{Kocic},
\cite{Milne}:
\bee
 x_{n+1} = \frac{a(n)x_n+b(n)}{c(n)x_n+d(n)},
 \label{wz2}
\eee
where the denominator cannot be zero.  We will first consider \eqref{rek3}.  In  our  case:
\bee
 a(n)=d(n)=(T-\gamma_n)^2+\left(\alpha-\frac12\right)\left(\beta-\frac12\right),
\eee
\bee
b(n)=-c(n)=(\alpha-\beta)(T-\gamma_n).
\eee
The customary way to solve the Ricccati  equation is to perform transformation of unknown variables in order to get second order linear  difference
equation.  We introduce new  variable
\bee
y(n)=H(n)f(n)
\eee
where
\bee
H(n)=\frac{C}{(T-\gamma_n)^2}, ~~~~C>1.
\eee
Then  the  relation \eqref{rek3}  becomes
\bee
y(n+1)=\frac{A(n)y(n)+B(n)}{C(n)y(n)+D(n)}
\eee
where
\bee
{   ~~}\begin{cases}
{   ~~} & A(n)=a(n)H(n+1)  \\
{   ~~}& B(n)=b(n)H(n)+H(n+1)\\
{   ~~} & C(n)=-b(n) \\
{   ~~} & D(n)=a(n)H(n)
    \end{cases}
\eee
Now we perform the second  transformation:
\bee
\frac{z(n+1)}{z(n)}=C(n)y(n)+D(n)
\label{wz4}
\eee
This   change of variables reduces our starting equation \eqref{rek3} to a  linear second order equation:
\bee
z(n+2)-P(n)z(n+1)-R(n)z(n)=0
\eee
where  we  have
\bee
{   ~~}\begin{cases}
{   ~~} & P(n)=D(n+1)+A(n)\displaystyle{\frac{C(n+1)}{C(n)}}\\
{   ~~} & R(n)=\left(B(n)C(n)-A(n)D(n)\right)\displaystyle{\frac{C(n+1)}{C(n)}}
    \end{cases}
\eee
We note that
\bee
\lim_{n\to \infty} P(n)=2C,  ~~~~\lim_{n\to \infty} R(n)=-C^2
\eee
so the characteristic equation for the linear  linear second order difference  equation is:
\bee
\lambda^2-2C\lambda +C^2=0
\eee
and the above equation has a double root  $\lambda=C$.  Going back to \eqref{wz4} we have;
\bee
-b(n)H(n)x(n)+a(n)H(n)=\frac{z(n+1)}{z(n)}
\label{wz8}
\eee
From Perron's  theorem  see  \cite[Chapt.17]{Milne}    we have
\bee
\limsup_{n \to \infty}|z(n)|^{1/n}=C
\eee
We also know that
\bee
\liminf_{n \to \infty}\frac{|z(n+1)|}{|z(n)|}\leq \liminf_{n \to \infty}|z(n)|^{1/n}\leq \limsup_{n \to \infty}|z(n)|^{1/n}\leq \limsup_{n \to \infty}\frac{|z(n+1)|}{|z(n)|}
\eee
so we have
\bee
C=\limsup_{n \to \infty}|z(n)|^{1/n}\leq  \limsup_{n \to \infty}\frac{|z(n+1)|}{|z(n)|},
\label{wz8}
\eee
\bee
\lim_{n \to \infty} a(n)H(n)=C.
\eee
From above and \eqref{wz8} we  conclude that
\bee
\lim_{n \to \infty }\frac{|x(n)|}{|T-\gamma_n|}=0.
\label{wz11}
\eee
We also note that
\bee
x(n+1)-x(n)=\frac{b(n)(1+(x(n)^2)}{a(n)-b(n)x(n)}.
\eee
From \eqref{wz11}  and above equation we obtain  that  for sufficiently  large  $n$ we have
\bee
x(n+1)-x(n)>0
\eee
so  for  large $n$  the  sequence $x(n)$  is increasing.  From \eqref{wz2} we see that  $x(n)$  cannot have real limit, because the
equation  fulfilled by the limit $x^\ast=\lim_{n \to \infty} x(n)$
\bee
 x^\ast  = \frac{a(n)x^\ast +b(n)}{-b(n)x^\ast + a(n)}
 \eee
has  no real  solution.  The conclusion is that the sequence $x(n)$  tends to infinity and  hence  we have also
\bee
\lim_{n \to \infty} f(n)=\infty
\eee
what we wanted to obtain.   Similar  calculations lead to the conclusion
\bee
\lim_{n \to \infty} g(n)=-\infty
\eee
and the Lemma 2  is  proved.

\end{document}